\documentclass{amsart}
\usepackage{amssymb}
\newtheorem{theorem}{Theorem}[section]
\newtheorem{lemma}[theorem]{Lemma}
\newtheorem{corollary}[theorem]{Corollary}
\newtheorem{proposition}[theorem]{Proposition}
\theoremstyle{definition}
\newtheorem{definition}[theorem]{Definition}
\newtheorem{example}[theorem]{Example}
\theoremstyle{remark}
\newtheorem{remark}[theorem]{Remark}
\def\htt{{\rm ht}}
\def\Max{{\rm Max}}
\def\Spec{{\rm Spec}}

\def\qf{{\rm qf}}
\def\Q{{\mathbb{Q}}}
\def\Z{{\mathbb{Z}}}

\begin{document}
%%%%%%%%%%%%%%%%%%%%%%%%%%%%%%%%%%%%%%%%%%%%%%%%%%%%%%%%%%%%
%%%%%%%%%%%%%%%%%%%%%%%%%%%%%%%%%%%%%%%%%%%%%%%%%%%%%%%%%%%%
\vspace*{1.5cm}
\begin{center}
\hrule\medskip
 {\tt Journal of Algebra 264 (2) (2003) 620-640 }
\medskip\hrule
\vspace{1cm}

 {\LARGE\bf Class Semigroups of Integral Domains} \vspace{1.5cm}

{\Large\bf S. Kabbaj $^{a,}$\footnote{ I would like to thank
Benedict H. Gross for helpful discussions. Research supported by
the Arab Fund for Economic and Social Development. Permanent
address: Department of Mathematics, P.O. Box 5046, KFUPM, Dhahran
31261, Saudi Arabia.} and A. Mimouni $^b$}
\vspace{1cm}

{\small $^a$ Department of Mathematics, Harvard University,
Cambridge, MA 02138, USA. \\ E-mail address:
kabbaj@math.harvard.edu \medskip

$^b$ Department of Mathematics, University of Fez, Fez 30000,
Morocco. \\E-mail address: a\_mimouni@hotmail.com}
\end{center}
\vspace{1.5cm}

\hrule\medskip
%%%%%%%%%%%%%%%%%%%%%%%%%%%%%%%%%%%%%%%%%%%%%%%%%%%%%%%%%%%%
%%%%%%%%%%%%%%%%%%%%%%%%%%%%%%%%%%%%%%%%%%%%%%%%%%%%%%%%%%%%
%%%Abstract%%%%%%%%%%%%%%%%%%%%%%%%%%%%%%%%%%%%%%%%%%%%%%%%%
%%%%%%%%%%%%%%%%%%%%%%%%%%%%%%%%%%%%%%%%%%%%%%%%%%%%%%%%%%%%
\noindent{\large\bf Abstract} \medskip

{\small This paper seeks ring-theoretic conditions of an integral
domain $R$ that reflect in the Clifford property or Boolean
property of its class semigroup ${\mathcal S}(R)$, that is, the
semigroup of the isomorphy classes of the nonzero (integral)
ideals of $R$ with the operation induced by multiplication.
Precisely, in Section 3, we characterize integrally closed domains
with Boolean class semigoup; in this case, ${\mathcal S}(R)$
identifies with the Boolean semigroup formed of all fractional
overrings of $R$. In Section 4, we investigate Noetherian-like
settings where the Clifford and Boolean properties of ${\mathcal
S}(R)$ coincide with (Lipman and Sally-Vasconcelos) stability
conditions; a main feature is that the Clifford property forces
$t-$locally Noetherian domains to be one-dimensional Noetherian
domains. Section 5 studies the transfer of the Clifford and
Boolean properties to various pullback constructions. Our results
lead to new families of integral domains with Clifford or Boolean
class semigroup, moving therefore beyond the contexts of
integrally closed domains or Noetherian domains.}
\medskip\hrule
%%%%%%%%%%%%%%%%%%%%%%%%%%%%%%%%%%%%%%%%%%%%%%%%%%%%%%%%%%%%
%%%%%%%%%%%%%%%%%%%%%%%%%%%%%%%%%%%%%%%%%%%%%%%%%%%%%%%%%%%%
%\subjclass[2000]{Primary 13C20, 13F05, 13F30, 13A15, 13B22,
%13G05; Secondary 20M14, 11R29}
%\keywords{Class semigroup, Clifford semigroup, Boolean semigroup,
%stable domain, Bezout domain, Pr\"ufer domain, Noetherian domain,
%Mori domain, pullbacks.}
%%%%%%%%%%%%%%%%%%%%%%%%%%%%%%%%%%%%%%%%%%%%%%%%%%%%%%%%%%%%
\newpage
%%%%%%%%%%%%%%%%%%%%%%%%%%%%%%%%%%%%%%%%%%%%%%%%%%%%%%%%%%%%
%%%Section1%%%%%%%%%%%%%%%%%%%%%%%%%%%%%%%%%%%%%%%%%%%%%%%%%
\begin{section}{Introduction}

Let $R$ be an integral domain. Following \cite{ZZ}, we define the
class semigroup of $R$, denoted ${\mathcal S}(R)$, to be the
(multiplicative Abelian) semigroup of nonzero fractional ideals
modulo its subsemigroup of nonzero principal ideals. The class
semigroup of $R$ contains, as subgroups, the class group Cl$(R)$
and, hence, the Picard group Pic$(R)$ of $R$.

In 1994, Zanardo and Zannier \cite{ZZ} proved that if $R$ is an integrally closed
domain and ${\mathcal S}(R)$ is a Clifford semigroup then $R$ is a Pr\"ufer
domain. The converse is not true since they showed that the ring of all
entire functions in the complex plane (which is Bezout) fails to have this
property. Their main result states that all orders in quadratic fields have
Clifford class semigroup. In 1996, Bazzoni and Salce \cite{BS1} investigated the structure
of the class semigroup for a valuation domain $V$, stating that ${\mathcal S}(V)$ is a
Clifford semigroup. In \cite{Ba1} and \cite{Ba2}, Bazzoni examined the case of Pr\"ufer
domains of finite character, showing that these have Clifford class semigroup, too. Recently,
she proved the converse in the case of integrally closed domains \cite{Ba3}.

This paper aims at investigating ring-theoretic properties
of an integral domain $R$ which reflect in the Clifford property or the
Boolean property of ${\mathcal S}(R)$. Precisely, in Section 3, our main
theorem asserts that ``{\em an integrally closed domain $R$ has Boolean class
semigroup if and only if $R$ is a strongly discrete Bezout domain of finite character
if and only if each nonzero ideal of $R$ is principal in its endomorphism ring}.'' One
may view this result as a satisfactory analogue of both \cite[Theorem 4.5]{Ba3} on the Clifford
property and \cite[Theorem 4.6]{O1} on stability. As a prelude to this, we characterize valuation
domains with Boolean class semigroup, stating that these are exactly the strongly discrete valuation
domains \cite{FP}. Section 4 studies Noetherian-like contexts. We prove that ``{\em if $R$ is a
$t-$locally Noetherian domain, then $R$ has Clifford (resp., Boolean) class semigroup if and only if $R$
is stable (resp., each nonzero ideal of $R$ is principal in its endomorphism ring)}.'' In particular, $t-$locally
Noetherian domains (such as Noetherian or strong Mori domains) with Clifford class semigroup turn out to be
one-dimensional Noetherian domains.
We also provide a characterization of Mori domains with Clifford or Boolean class semigroup that links them
to stability, specifically,  ``{\em a Mori domain $R$ is stable (resp., each nonzero ideal of $R$ is principal
in its endomorphism ring) if and only if $R$ is
a one-dimensional Clifford (resp., Boole) regular domain and the complete integral closure of $R$ is Mori}.''
Section 5 treats the possible transfer of the Clifford and Boolean properties to pullbacks. New families of
domains with Clifford or Boolean class semigroup stem from our results. Throughout, examples are provided to
illustrate the scopes and limits of the results.

For the convenience of the reader, we summarize in the following two diagrams the relations between the main
classes of domains involved in this paper (where ``+ IC'' means that the implication requires the integrally
closed hypothesis):\bigskip

%%DIAGRAM1%%%%%%%%%%%%%%%%%%%%%%%%%%
%%%%%%%%%%%%%%%%%%%%%%%%%%%%%%%%%%%
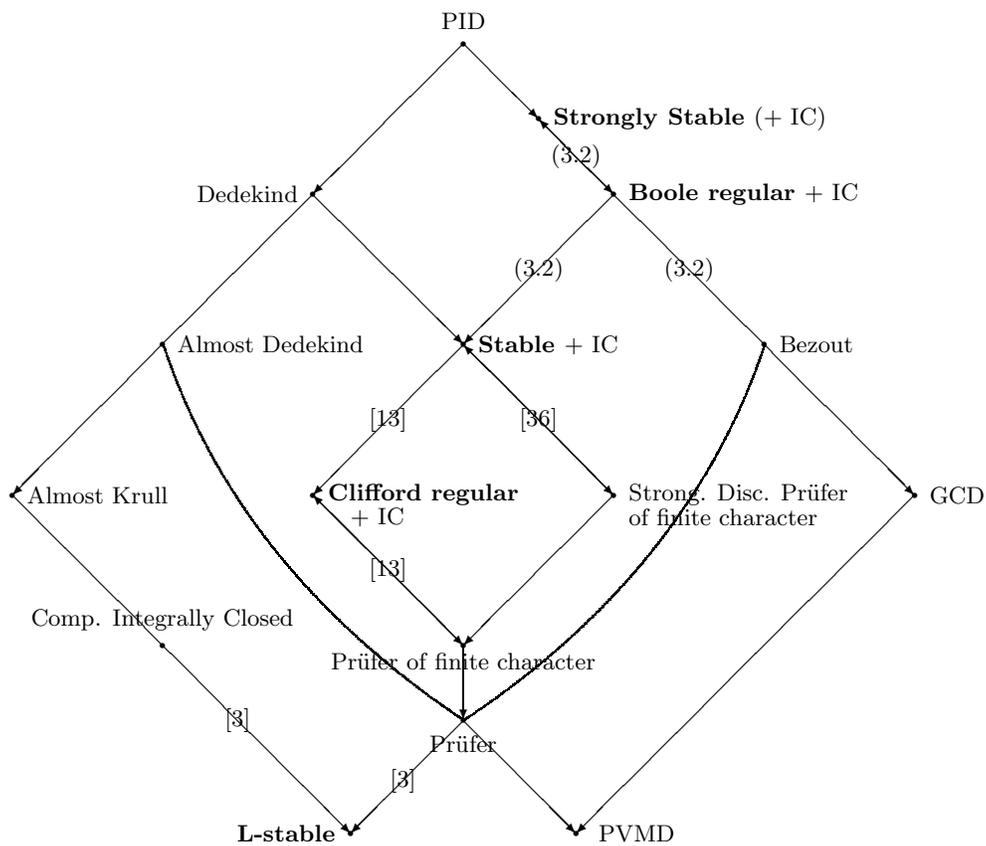
\[\setlength{\unitlength}{1mm}
\begin{picture}(120,110)(-15,-20)
%%%%%%%%%%%%
\put(40,90){\vector(1,-1){10}}
\put(50,80){\vector(1,-1){10}}
\put(40,90){\vector(-1,-1){20}}
\put(60,70){\vector(1,-1){40}}
\put(60,70){\vector(-1,-1){20}}
\put(60,70){\vector(-1,1){10}}
\put(40,50){\vector(1,-1){20}}
\put(40,50){\vector(-1,-1){20}}
\put(20,70){\vector(1,-1){20}}
\put(20,70){\vector(-1,-1){40}}
\put(20,30){\vector(1,-1){20}}
\put(60,30){\vector(-1,-1){20}}
\put(60,30){\vector(-1,1){20}}
\put(-20,30){\vector(1,-1){45}}
\put(100,30){\vector(-1,-1){45}}
\put(40,10){\vector(-1,1){20}}
\put(40,10){\vector(0,-1){10}}
\put(40,0){\vector(1,-1){15}}
\put(40,0){\vector(-1,-1){15}}
\qbezier(0,50)(10,20)(40,0)
\qbezier(80,50)(70,20)(40,0)
%%%%%%%%%%%%%%%%%
\put(40,90){\circle*{.7}}
\put(40,92){\makebox(0,0)[b]{\small PID}}
\put(50,80){\circle*{.7}}
\put(52,80){\makebox(0,0)[l]{\small\bf Strongly Stable {\rm (+ IC)}}}
\put(60,70){\circle*{.7}}
\put(62,70){\makebox(0,0)[l]{\small \bf Boole regular {\rm + IC}}}
\put(20,70){\circle*{.7}}
\put(18,70){\makebox(0,0)[r]{\small Dedekind}}
\put(40,50){\circle*{.7}}
\put(42,50){\makebox(0,0)[l]{\small\bf Stable {\rm + IC}}}
\put(0,50){\circle*{.7}}
\put(2,50){\makebox(0,0)[l]{\small Almost Dedekind}}
\put(20,30){\circle*{.7}}
\put(22,30){\makebox(0,0)[l]{\small \bf Clifford regular}}
\put(25,27){\makebox(0,0)[l]{\small + IC}}
\put(60,30){\circle*{.7}}
\put(62,30){\makebox(0,0)[l]{\small Strong. Disc. Pr\"ufer}}
\put(62,27){\makebox(0,0)[l]{\small of finite character}}
\put(100,30){\circle*{.7}}
\put(102,30){\makebox(0,0)[l]{\small GCD}}
\put(80,50){\circle*{.7}}
\put(82,50){\makebox(0,0)[l]{\small Bezout}}
\put(-20,30){\circle*{.7}}
\put(-18,30){\makebox(0,0)[l]{\small Almost Krull}}
\put(40,10){\circle*{.7}}
\put(40,9){\makebox(0,0)[t]{\small Pr\"ufer of finite character}}
\put(40,0){\circle*{.7}}
\put(40,-2){\makebox(0,0)[t]{\small Pr\"ufer}}
\put(0,10){\circle*{.7}}
\put(0,12){\makebox(0,0)[b]{\small Comp. Integrally Closed}}
\put(25,-15){\circle*{.7}}
\put(23,-15){\makebox(0,0)[r]{\small\bf L-stable}}
\put(55,-15){\circle*{.7}}
\put(58,-15){\makebox(0,0)[l]{\small PVMD}}
%%%%%%%%%%%%%%%%%%%%%
\put(50,40){\makebox(0,0){\small \cite{O1}}}
\put(30,20){\makebox(0,0){\small \cite{Ba3}}}
\put(30,40){\makebox(0,0){\small \cite{Ba3}}}
\put(32,-8){\makebox(0,0){\small \cite{AHP}}}
\put(10,0){\makebox(0,0){\small \cite{AHP}}}
\put(70,60){\makebox(0,0){\small (\ref{3.1})}}
\put(50,60){\makebox(0,0){\small (\ref{3.1})}}
\put(55,75){\makebox(0,0){\small (\ref{3.1})}}
%%%%%%%%%%%%%%%%%%%%%
\end{picture}\]\bigskip
%%%%%%%%%%%%%%%%%%%%%

%%DIAGRAM2%%%%%%%%%%%%%%%%%%%%%%%%%%
%%%%%%%%%%%%%%%%%%%%%%%%%%%%%%%%%%%
\[\setlength{\unitlength}{1mm}
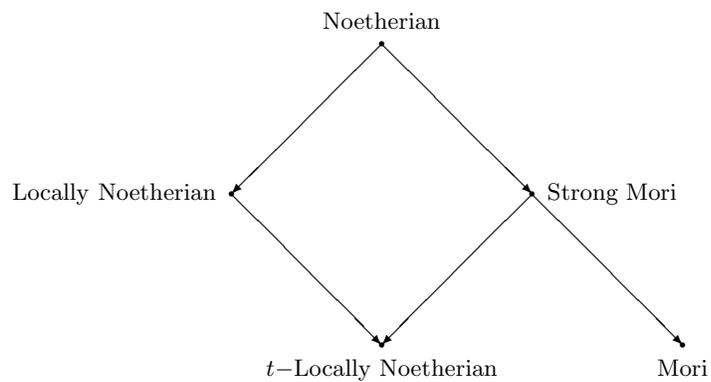
\begin{picture}(60,50)(18,-80)
%%%%%%%%%%%%
\put(40,-30){\vector(1,-1){20}}
\put(40,-30){\vector(-1,-1){20}}
\put(20,-50){\vector(1,-1){20}}
\put(60,-50){\vector(-1,-1){20}}
\put(60,-50){\vector(1,-1){20}}
%\qbezier(0,50)(-55,30)(20,-50)
%\qbezier(20,70)(-90,0)(40,-30)
%%%%%%%%%%%%%%%%%%%%%
\put(40,-30){\circle*{.7}}
\put(40,-28){\makebox(0,0)[b]{\small Noetherian}}
\put(20,-50){\circle*{.7}}
\put(18,-50){\makebox(0,0)[r]{\small Locally Noetherian}}
\put(60,-50){\circle*{.7}}
\put(62,-50){\makebox(0,0)[l]{\small Strong Mori}}
\put(40,-70){\circle*{.7}}
\put(40,-72){\makebox(0,0)[t]{\small  $t-$Locally Noetherian}}
\put(80,-70){\circle*{.7}}
\put(80,-72){\makebox(0,0)[t]{\small  Mori}}
\end{picture}\]
%%%%%%%%%%%%%%%%%%%%%%%%%%%%%%%%%%%%%%%%%%%%%%%%%%%%%%%%%
\end{section}
%%%%%%%%%%%%%%%%%%%%%%%%%%%%%%%%%%%%%%%%%%%%%%%%%%%%%%%%%
\newpage
%%%%%%%%%%%%%%%%%%%%%%%%%%%%%%%%%%%%%%%%%%%%%%%%%%%%%%%%%
%%%SECTION2%%%%%%%%%%%%%%%%%%%%%%%%%%%%%%%%%%%%%%%%%%%%%
\section{Preliminaries}

Let us first recall  the following definitions. A commutative
semigroup $S$ is said to be a {\em Clifford semigroup} if every
element $x$ of $S$ is (von Neumann) regular, i.e., there exists
$a\in S$ such that $x^{2}a=x$; and $S$ is said to be {\em Boolean}
if for each $x\in S$, $x=x^{2}$ (cf. \cite{Ho}). The importance of
a Clifford semigroup $S$ resides in its ability to stand as a
disjoint union of subgroups $G_e$, where $e$ ranges over the set
of idempotent elements of $S$, and $G_e$ is the largest subgroup
of $S$ with identity equal to $e$. Often, the $G_e$'s are called
the constituent groups of $S$. Clearly, a semigroup $S$ is Boolean
if and only if the constituent groups of $S$ are all trivial.

As in \cite{Ba3}, we say that a domain $R$ is {\em Clifford regular} if the class semigroup
${\mathcal S}(R)$ of $R$ is a Clifford semigroup. By analogy with this,  we say that a domain
$R$ is {\em Boole regular} if the class semigroup ${\mathcal S}(R)$ of $R$ is a Boolean semigroup.
At this point, recall Bazzoni's recent result \cite[Theorem 4.5]{Ba3}: {\em an integrally closed
domain $R$ is Clifford regular if and only if $R$ is a Pr\"ufer domain of finite character}
(i.e., each nonzero ideal is contained only in finitely many maximal ideals).

An ideal of an integral domain $R$ is said to be {\em L-stable}
(here L stands for Lipman) if $R^{I} =\bigcup (I^{n}:I^{n}) =
(I:I)$, and $R$ is called an {\em L-stable domain} if every
nonzero ideal of $R$ is L-stable \cite{AHP}. Lipman \cite{Li}
introduced the notion of stability in the specific setting of
one-dimensional commutative semi-local Noetherian  rings (to give
a characterization of Arf rings). In Lipman's context, {\em an
integral domain $R$  is L-stable if and only if $R$ is Boole
regular} (cf. \cite[Lemma 1.11]{Li}).

An ideal $I$ of an integral domain $R$ is said to be {\em stable}
if $I$ is invertible in $(I:I)$, and $R$ is called a {\em stable
domain} provided each nonzero ideal of $R$ is stable \cite{AHP}.
Sally and Vasconcelos \cite{SV} used this concept to settle Bass'
conjecture on one-dimensional Noetherian rings with finite
integral closure. Recall that a stable domain is L-stable
\cite[Lemma 2.1]{AHP}. For recent developments on stability (in
settings different than originally considered), we refer the
reader to \cite{AHP,Ba3,O1,O2,O3}.  Of particular relevance to our
study is Olberding's result \cite[Theorem 4.6]{O1} stating that
{\em an integrally closed domain $R$ is stable if and only if $R$
is a strongly discrete Pr\"ufer domain of finite character}.

Throughout, all rings considered are integral domains. We shall use ${\overline I}$ to denote   the isomorphy class of an ideal I.
\bigskip

We often will be appealing to the next results without explicit mention.

%%%2.1%%%%%%%%%%%%%%%%%%%%%%%%%%%%%%%%%%%%%%%%
\begin{lemma} \label{2.1}
1) Let $I$ be an ideal of an integral domain $R$. ${\overline I}$ is a regular element of ${\mathcal S}(R)$ if and only if $I =I^{2}(I:I^{2})$ {\rm (\cite[Lemma 1.1]{Ba1})}. \\
2) A stable domain is Clifford regular {\rm (\cite[Proposition 2.2]{Ba3})}.\\
3) A stable domain has finite character {\rm (\cite[Theorem 3.3]{O3})}. \\
4) An integrally closed stable domain is Pr\"ufer {\rm (\cite[Lemma F]{ES})}.
\end{lemma}

The next lemma establishes the transfer of the Clifford and Boolean properties to two types of overrings.

%%%2.2%%%old1.4%%%%%%%%%%%%%%%%%%%%%%%%%%%%%%%
\begin{lemma} \label{2.3}
Let $R$ be an integral domain and $B$ an overring
of $R$. Assume that one of the following two assumptions holds:\\
a) $B$ is a flat extension of $R$,\\
b) The conductor $(R:B)$ is nonzero.\\
If $R$ is a Clifford (resp., Boole) regular domain, then so is
$B$.
\end{lemma}
%%%%%%%%%%%%%%%%%%%%%%
 \begin{proof}[\bf Proof.] a) Let $J$ be an ideal of $B$. It suffices to
show that $J\subseteq J^{2}(J:J^{2})$. Let $I:=J\cap R$. By
\cite[Proposition 1.2 (ii)]{P}, $J=IB$. For each $x\in (I:I^{2})$,
$xI^{2}\subseteq I$ implies that $xI^{2}B\subseteq IB$. Hence
$xJ^{2} = x(IB)^{2} = xI^{2}B \subseteq IB = J$. So $x\in
(J:J^{2})$ and hence $(I:I^{2})\subseteq (J:J^{2})$. Therefore $I
= I^{2}(I:I^{2})\subseteq J^{2}(J:J^{2})$. So that $J\subseteq
J^{2}(J:J^{2})$.\\
\noindent b) Assume that $(R:B)\not =0$. Let $c\in
(R:B)\setminus{0}$, $J$ an ideal of $B$, and  $I=cJ$. Clearly, $I$
is an ideal of $R$ with $I^{2}(I:I^{2})=cJ^{2}(J:J^{2})$. Hence
$cJ = I = I^{2}(I:I^{2}) = cJ^{2}(J:J^{2})$. It follows that $J =
J^{2}(J:J^{2})$ and hence ${\overline J}$ is regular in ${\mathcal
S}(B)$. Consequently, $B$ is Clifford regular. Now assume $R$ is
Boole regular. Here it suffices to notice that if $I^{2}=qI$, then
$J^{2}=qJ$. \end{proof}

 Our next result, Proposition~\ref{2.4}, will play a central role in the development of Sections 3 and 4. It generalizes Zanardo-Zannier's theorem mentioned in the introduction.

%%%2.3%%%old1.6%%%%%%%%%%%%%%%%%%%%%%%%%%%%
\begin{proposition} \label{2.4}
Let $R$ be an integral domain. If $R$ is a
Clifford (resp., Boole) regular domain, then ${\overline R}$ is a Pr\"ufer (resp.,
Bezout) domain, where ${\overline R}$ denotes the integral closure of $R$.
\end{proposition}
%%%%%%%%%%%%%%%%%%%%%%
 \begin{proof}[\bf Proof.] The Clifford statement is handled by \cite[Proposition 2.1]{R}
 and \cite[Proposition 2.3]{Ba3}. Next assume that $R$ is a Boole regular domain.  By the first part, ${\overline R}$ is a Pr\"ufer domain. Let $J$ be a
finitely generated ideal of ${\overline R}$. Write $J=\sum_{i=1}^{i=r}a_{i}{\overline R}$. Let
$T:=R[a_{1}, \dots, a_{r}]$ and $I:=\sum_{i=1}^{i=r}a_{i}T$. Since $T$ is a
finitely generated $R$-module, then $(R:T)\not =0$. By Lemma~\ref{2.3}, ${\mathcal S}(T)$
is Boolean. So there is $0\not =c\in K$ such that $I^{2}=cI$. Since $I{\overline R}=J$,
then $J^{2}=cJ$. Hence $(J:J^{2})=(J:cJ)=c^{-1}(J:J)$. Since $J$ is
invertible in ${\overline R}$, then $(J:J)={\overline R}$, hence
$c^{-1}{\overline R}=c^{-1}(J:J)=(J:J^{2})=((J:J):J)=({\overline R}:J)$, whence
$c^{-1}J=J({\overline R}:J)={\overline R}$. So $J=c{\overline R}$ and thus ${\overline R}$ is a Bezout domain. \end{proof}

Our first corollary characterizes almost Krull domains with Clifford or Boolean class semigroup. Notice that our elementary proof of this result does not appeal to \cite[Theorem 4.5]{Ba3}, rather it draws on basic properties of almost Krull domains.

%%%2.4%%%old1.8%%%%%%%%%%%%%%%%%%%%%%%%%%%%
\begin{corollary} \label{2.6}
Let $R$ be an integral domain. Then  $R$ is almost Krull and Clifford (resp., Boole) regular if and only if $R$ is Dedekind (resp., a PID).
\end{corollary}
%%%%%%%%%%%%%%%%%%%%%
\begin{proof}[\bf Proof.] We just need to prove the ``only if" assertion. Clearly, for any maximal ideal $M$ of $R$,
$R_{M}$ inherits the Clifford property from $R$. Hence, by Proposition~\ref{2.4} $R$ is an almost Dedekind domain. Suppose that
there exists a nonzero ideal $I$ of $R$ which is not invertible, i.e.,
$II^{-1}\subsetneqq R$.
Let $J := II^{-1}$. Then $J$ is a proper trace ideal of
$R$, hence $J^{-1} = (J:J) =R$ (since $R$ is
completely integrally closed), whence $ (J:J^{2}) = ((J:J):J)
= (R:J) = J^{-1} = R$. So $J=J^{2}(J:J^{2}) = J^{2}$ (since  ${\overline J}$ is
regular in ${\mathcal S}(R)$). It follows that $J=J^{n}$, for each
$n\geq 1$. Since $R$ is almost Dedekind, $J = \bigcap_{n\geq 1}(J^{n}) = (0)$, the desired contradiction.

The Boolean statement follows from the Clifford statement and Proposition~\ref{2.4}, completing the proof. \end{proof}

A brief discussion at the end of Section 3 envisages a possible widening of the scope of Corollary~\ref{2.6} to completely integrally closed domains.

%%%2.5%%%old1.9%%%%%%%%%%%%%%%%%%%%%%%%%%%%
\begin{corollary} \label{2.7}
Let $R$ be an integral domain and $X$ an indeterminate
over $R$. The following statements are equivalent:\\
i) $R$ is a field;\\
ii) $R[X]$ is Boole regular;\\
iii) $R[X]$ is Clifford regular.
\end{corollary}

%%%%%%%%%%%%%%%%%%%%%%%%%%%%%%%%%%%%%%%%%%%%%%%%%%%%%%%
%%%SECTION3%%%%%%%%%%%%%%%%%%%%%%%%%%%%%%%%%%%%%%%%%%%%%
%%%%%%%%%%%%%%%%%%%%%%%%%%%%%%%%%%%%%%%%%%%%%%%%%%%%%%%%
\section{Boole regular domains}

Clearly, a PID is Boole regular (see definition in Section 2) and
a Boole regular domain is Clifford regular. Our purpose in this
section is to characterize Boole regularity for integrally closed
domains. Recall that the study of Clifford regularity -in the
integrally closed context- was initiated in \cite{Ba1,Ba2} and
recently achieved in \cite{Ba3}.

As a prelude, we characterize valuation domains with Boolean class
semigroup, stating that these are exactly the strongly discrete
valuation domains \cite{FP}. An integral domain is {\em strongly
discrete} if it has no nonzero idempotent prime ideals. A stable
domain trivially is strongly discrete.\bigskip

We shall first find a natural stability condition that best suits the Boolean context.
It can be termed as follows:

%%%3.1%%%%%%%%%%%%%%%%%%%%%%%%%%%%%%%
\begin{definition} \label{3.0}
An integral domain $R$ is called a {\em strongly stable domain} if
each nonzero ideal of $R$ is principal in its endomorphism ring
$(I:I)$.
\end{definition}

Next, we announce the main result of this section. First note that for any integral domain $R$,
the set  ${\mathcal F}_{OV}(R)$ of fractional {\em overrings} of $R$ is a Boolean semigroup with identity equal to $R$.

%%%3.2%%%old2.5%%%%%%%%%%%%%%%%%%%%%%%%%%%%
\begin{theorem} \label{3.1}
Let $R$ be an integrally closed domain. The following statements are equivalent:\\
i) $R$ is a Boole regular domain;\\
ii) $R$ is a strongly discrete Bezout domain of finite character;\\
iii) $R$ is a strongly stable domain.

Moreover, when any one condition holds, ${\mathcal S}(R)={\mathcal F}_{OV}(R)$, where $\overline{T}$ is identified with $T$ for each fractional overring $T$ of $R$.
\end{theorem}
%%%%%%%%%%%%%%%%%

The proof involves some preliminary results of independent interest.

%%%3.3%%%old1.3%%%%%%%%%%%%%%%%%%%%%%%%%%%%
\begin{lemma} \label{3.2}
Let $R$ be an integral domain.  The following statements are equivalent:\\
i) $R$ is a stable Boole regular domain;\\
ii) $R$ is a strongly stable domain.
\end{lemma}
%%%%%%%%%%%%%%%%%%%%
\begin{proof}[\bf Proof.] $i)\Longrightarrow ii)$ Let $I$ be a nonzero ideal of
$R$. Since ${\mathcal S}(R)$ is Boolean, then $I^{2}=cI$ for some $0\not =c\in K$. So
$(I:I^{2})=(I:cI)=c^{-1}(I:I)$. Since $R$ is stable, then
$I(I:I^{2})=(I:I)$. Hence $c^{-1}I=I(I:I^{2})=(I:I)$ and therefore
$I=c(I:I)$.\\
$ii)\Longrightarrow i)$ Clearly, $R$ is stable. Further, let $I$ be a nonzero ideal of $R$. If $I=c(I:I)$,
then $I^{2}=cI$, as desired. \end{proof}

%%%3.4%%%old2.2%%%%%%%%%%%%%%%%%%%%%%%%%%%%
\begin{lemma} \label{3.3}
Let $R$ be an integrally closed domain. The following statements are equivalent:\\
i) $R$ is a strongly discrete Clifford regular domain;\\
ii)  $R$ is a stable domain.
\end{lemma}
%%%%%%%%%%%%%%%%%%%%%
\begin{proof}[\bf Proof.] By \cite{O1} we need only prove $(i)\Longrightarrow (ii)$. This follows from a combination of \cite[Theorem 4.5]{Ba3} and \cite[Theorem 4.6]{O1}; however, we offer the following different elementary proof (which draws on the basic fact that the maximal ideal of a strongly discrete valuation domain is principal \cite[Lemma 2.1]{FP}). Assume that (i) holds. By Proposition~\ref{2.4}, $R$ is a strongly discrete Pr\"ufer domain. Let $I$ be a nonzero ideal of $R$, $T:=(I:I)$, and
$J:=I(T:I)$. Since $\overline{I}$ is regular in ${\mathcal S}(R)$, then $I=IJ$ and
$J^{2}=J$ \cite[Proposition 2.1(1)]{Ba1}. Suppose that $J\subsetneqq T$. Let $Q$ be
a minimal prime ideal of $T$ over $J$ and $q=Q\cap R$. Then $T_{Q}=R_{q}$ is a strongly discrete valuation domain and hence
$QT_{Q}=aT_{Q}$ for some $0\not =a\in Q$. Since $Q$ is minimal over $J$,
then $JT_{Q}$ is $QT_{Q}$-primary. So $JT_{Q}=(QT_{Q})^{r}$, for some
integer
$r$. Since $J=J^{2}$, then $a^{r}T_{Q}=a^{2r}T_{Q}$, the desired contradiction. Therefore $J=T$ and hence $R$ is stable. \end{proof}

Recall that Bazzoni and Salce \cite{BS1} proved that valuation domains have
 always Clifford class semi group; next we characterize those among them with Boolean class semigoup.

%%%3.5%%%old2.4%%%%%%%%%%%%%%%%%%%%%%%%%%%%
\begin{lemma} \label{3.5}
 Let $V$ be a valuation domain. The following
assertions are equivalent:\\
i) $V$ is a Boole regular domain;\\
ii) $V_{P}$ is a divisorial domain, for each nonzero prime ideal
$P$ of $R$;\\
iii) $V$ is a stable domain;\\
iv) $V$ is a strongly discrete valuation domain.
\end{lemma}
%%%%%%%%%%%%%%%%%%%%%
\begin{proof}[\bf Proof.] $i)\Longrightarrow ii)$ { \sl Claim: If ${\mathcal S}(V)$
is Boolean, then $V$ is a divisorial domain}. Indeed, let $I$ be a nonzero
ideal of
$V$ and $Z(V, I)$ the set of zerodivisors of $R$ modulo $I$. Then $Z(V,
I):=P$
is a prime ideal of $V$ and $(I:I)=V_{P}$. Since ${\mathcal S}(V)$ is Boolean,
then there is $0\not =c\in K$ such that $I^{2}=cI$. Two cases are possible.
{\sl {Case 1:}} $I(V_{P}:I)=V_{P}$. Then $I=aV_{P}$, for some
nonzero $a\in I$. So $(V:I)=(V:aV_{P})= a^{-1}(V:V_{P})=a^{-1}P$. Hence
$I_{v}=(V:(V:I))=(V:a^{-1}P)=a(V:P)$. Now,
if $P$ is not a maximal ideal of $V$, then $(V:P)=(P:P)=V_{P}$; hence
$I_{v}=a(V:P)=aV_{P}=I$. So $I$ is divisorial. If $P$ is maximal in $V$,
then $I=aV$.
Here too, $I$ is divisorial.\\
{\sl {Case 2:}} $I(V_{P}:I)\subsetneqq V_{P}$. Since $V_{P}$ is a
$TP$-domain \cite{FHP1}, then there is a prime ideal $Q$ of $V$ with $Q\subseteq P$ such
that $I(V_{P}:I)=QV_{P}$. On the other hand, $I^{2}=cI$ yields
$(V_{P}:I)=(I:I^{2})=(I:cI)=c^{-1}V_{P}$. So that
$QV_{P}=I(V_{P}:I)=Ic^{-1}V_{P}=c^{-1}I$, whence $I=cQV_{P}$. So
$V_{P}=(I:I)=(cQV_{P}:cQV_{P})=(QV_{P}:QV_{P})=V_{Q}$. It follows that $P=Q$
and
$I=cQV_{P}=cPV_{P}=cP$. Since $I^{2}=cI$, then $P=P^{2}$. Now $P$
is a trace ideal of $V$. Then $(V:P)=(P:P)=V_{P}$. So
$(V:I)=(V:cP)=c^{-1}(V:P)=c^{-1}V_{P}$. Therefore
$I_{v}=(V:c^{-1}V_{P})=c(V:V_{P})=cP=I$ and hence $I$ is divisorial. Consequently, $V$ is divisorial, completing the proof of our claim.

Now, let $P$ be any nonzero prime ideal of
$V$. By Lemma~\ref{2.3}, $V_P$ inherits the Boolean property from $V$. By the above claim, $V_{P}$ is divisorial, as desired.\\
$ii)\Longrightarrow iii)$ Let $P$ be a prime ideal of $V$. By \cite[Lemma 5.2]{H}, $P=PV_{P}=aV_{P}$, for some $a\in P$. By \cite[Proposition 2.10]{AHP},
$V$ is stable.\\
$iii)\Longrightarrow i)$ Let $I$ be a nonzero ideal of $V$ and
$P:=Z(V, I)$. By (iii), $I$ is invertible in $(I:I)=V_{P}$. Hence
$I=aV_{P}$,
for some $a\in I$. So $I^{2}=aI$. Hence ${\mathcal S}(V)$ is Boolean. \\
$(iii)\Longleftrightarrow (iv)$ is handled by \cite[Proposition 2.10]{AHP}.
\end{proof}

Notice that Lemma~\ref{3.5} gives rise to a large class of Boole regular domains that are not PIDs. Indeed, any strongly discrete valuation domain of dimension $\geq 2$ does (e.g., $k[X]_{(X)}+Yk(X)[[Y]]$, where $k$ is a field and $X,Y$ are indeterminates over $k$  \cite{FP}).

%%%3.6%%%new%%%%%%%%%%%%%%%%%%%%%%%%%%%%
\begin{lemma} \label{3.6}
An integrally closed domain $R$ is locally Boole regular if and only if $R$ is a strongly discrete Pr\"ufer domain.
\end{lemma}
%%%%%%%%%%%%%%%%%%%%%%%%%%%%%
\begin{proof}[\bf Proof.] Combine  Proposition~\ref{2.4} and Lemma~\ref{3.5}. \end{proof}

%%%3.7%%%new%%%%%%%%%%%%%%%%%%%%%%%%%%%%
\begin{lemma} \label{3.6(1)}
An integrally closed domain $R$ is Boole regular if and only if $R$ is a stable Bezout domain.
\end{lemma}
%%%%%%%%%%%%%%%%%%%%%%%%%%%%%
\begin{proof}[\bf Proof.] Assume $R$ is Boole regular. By Proposition~\ref{2.4}, $R$ is Bezout. Further, a combination of Lemma~\ref{2.3} and Lemma~\ref{3.6} ensures that $R$ is a strongly discrete Pr\"ufer domain. It turns out that $R$ is a strongly discrete Clifford domain, hence it is stable by Lemma~\ref{3.3}. Conversely, Let $I$ be an ideal of $R$. Then $T:=(I:I)$ is a Bezout domain. Further, $I$ is invertible in $T$, so it is principal in $T$ to complete the proof.
\end{proof}

%%%%%%%%%%%%%%%%%%%%%%%%%%%%%%%%%%
\begin{proof}[\bf Proof of Theorem~\ref{3.1}.] $(i)\Longrightarrow (ii)$ Follows from Lemma~\ref{3.6(1)} along with the facts that a stable domain is necessarily strongly discrete and has finite character.\\
$(ii)\Longrightarrow (i)$ Follows from \cite[Theorem 4.6]{O1} (and Lemma~\ref{3.6(1)}); however, we offer the following direct proof which draws on Bazzoni's study of the groups associated to idempotents in the class semigroup. Next, assume that $R$ is a strongly discrete Bezout domain of finite character. Then ${\mathcal S}(R)=\bigvee G_{\overline J}$, where $\overline{J}$ ranges over the set of idempotent elements of ${\mathcal S}(R)$.  By \cite[Theorem 3.1]{Ba2}, an element
$\overline{J}$ of ${\mathcal S}(R)$ is idempotent if and only if there exists a unique nonzero idempotent fractional ideal $L$ of $R$ such that $J\cong L$ and $L$ satisfies one of the following two conditions:\\
(1) $L=T$, where $T$ is a fractional overring of $R$, or\\
(2) $L=P_{1}P_{2}...P_{n}T$, where each $P_{i}$ is a nonzero idempotent
prime ideal of $R$ and $T$ is a fractional overring of $R$.\\
Since $R$ is strongly discrete, then there is no nonzero idempotent prime ideals. This rules out the $L$'s issued from the second condition. Further, by \cite[Proposition 2.2]{Ba4}, the group $G_{\overline{T}}$ associated to $\overline{T}$ coincides with the class group Cl$(T)$ for each fractional overring $T$ of $R$. Since $R$ is Bezout, then each overring $T$ of $R$ is Bezout and therefore Cl$(T)$ is trivial. Hence the constituent groups of ${\mathcal S}(R)$ are all trivial, whence ${\mathcal S}(R)$ is Boolean, as desired.\\
$(i)\Longleftrightarrow (iii)$ is handled by Lemma~\ref{3.2} and Lemma~\ref{3.6(1)}.

Finally, assume that $(i)-(iii)$ hold. Clearly, ${\mathcal S}(R)=\{\overline{T}: T\in {\mathcal F}_{OV}(R)\}$ by \cite[Theorem 3.1]{Ba2} mentioned above. Moreover, due to the uniqueness required by this theorem, one can identify $\overline{T}$ with $T$ for each $T\in {\mathcal F}_{OV}(R)$, leading therefore to the identification of ${\mathcal S}(R)$ with  the Boolean semigroup ${\mathcal F}_{OV}(R)$, completing the proof of the theorem. \end{proof}

%%%3.8%%%new%%%%%%%%%%%%%%%%%%%%%%%%%%%%
\begin{example} \label{3.8}
In \cite[Construction 1]{Lo}, Loper shaped an example of a generalized Dedekind domain (hence a strongly discrete Pr\"ufer domain \cite{FHP2}) which is not Bezout. Further, (one can easily check that) it has finite character. Hence it is stable \cite{O1} but not Boole regular (Theorem~\ref{3.1}). It follows that Theorem~\ref{3.1} does not extend to strongly discrete Pr\"ufer domains of finite character (equivalently, integrally closed stable domains).
\end{example}

%%%3.9%%%new%%%%%%%%%%%%%%%%%%%%%%%%%%%%
\begin{remark} \label{3.7} Theorem~\ref{3.1} and its satellite lemmas yield immediate consequences:\\
1) Unlike Clifford regularity, Boole regularity is not a local property for the class of integrally closed domains of finite character.\\
2) If $R$ is an integrally closed domain that is Boole regular (equivalently, strongly stable), then so is any overring of $R$.\\
3) Stability and strong stability do not coincide in general (e.g., Dedekind domains that are not PIDs). They do however in integrally closed semilocal contexts (see Corollary~\ref{3.9}).\\
4) Unlike stability, strong stability is not a local property for the class of domains of finite character.\\
5)  If $R$ is a strongly stable domain, then so is its integral closure ${\overline R}$.
\end{remark}

Moreover, a Bezout domain of finite character need not be Boole regular (e.g., valuation domains with nonzero idempotent prime ideals). Consequently, in view of the above discussion, Theorem~\ref{3.1} may stand as a satisfactory analogue of both \cite[Theorem 4.5]{Ba3} and \cite[Theorem 4.6]{O1} for Boole regularity and strong stability, respectively. \bigskip

In the semilocal context where ``Pr\"ufer'' elevates to ``Bezout'', most of the notions in play collapse, as shown by the next result.

%%%3.10%%%%%%%%%%%%%%%%%%%%%%%%%%%%%%%
\begin{corollary} \label{3.9}
Let $R$ be an integrally closed semilocal domain. The following statements are equivalent:\\
i) $R$ is a strongly stable domain;\\
ii) $R$ is a Boole regular domain;\\
iii) $R$ is a stable domain;\\
iv) $R$ is a strongly discrete Clifford regular domain;\\
v) $R$ is a strongly discrete Pr\"ufer domain.
\end{corollary}
%%%%%%%%%%%%%%%%%

It is worth noticing that from Corollary~\ref{3.9} stems a large family of  examples of integrally closed Boole regular domains that are neither PIDs nor strongly discrete valuation domains (e.g., semilocal strongly discrete Pr\"ufer domains  of dimension $\geq 2$). Recall that the class of strongly discrete Pr\"ufer domains of finite character properly contains the class of integrally closed Boole regular domains. \bigskip

We close this section with a brief discussion of the completely integrally closed case. Indeed, by Theorem~\ref{3.1}, {\em a  completely integrally closed domain is Boole regular if and only if it is a PID}. This extends the Boolean statement of Corollary~\ref{2.6}. However, a one-dimensional completely integrally closed Clifford regular domain (e.g.,  a non-discrete rank-one valuation domain) need not be Dedekind. Compare to the Clifford statement of Corollary~\ref{2.6} as well as to the known fact that a one-dimensional integrally closed  stable domain is Dedekind.

%%%SECTION4%%%%%%%%%%%%%%%%%%%%%%%%%%%%%%%%%%%%%%%%%%%%%
%%%%%%%%%%%%%%%%%%%%%%%%%%%%%%%%%%%%%%%%%%%%%%%%%%%%%%%%
\section{Noetherian-like settings}

This section investigates the class semigroup for two large classes of Noetherian-like
domains, that is, $t-$locally Noetherian domains and Mori domains. Precisely, we study
conditions under which stability and strong stability characterize Clifford regularity
and Boole regularity, respectively. A main feature of our first theorem is that Clifford
regularity forces the Noetherianity of $t-$locally Noetherian domains. However, the second
main theorem (on Mori domains) may allow one, a priori, to move beyond the context of Noetherian
domains. Unfortunately, we are not able to shape an example that supports this claim. (See the
brief discussion at the end of this section.)

In order to provide some background for the present section, we
review  some terminology related to star-operations \cite{G}. Let
$R$ be an integral domain. For a nonzero fractional ideal $I$ of
$R$, set $I_v:=(I^{-1})^{-1}$; $I_t:=\bigcup J_v$ where $J$ ranges
over the set of finitely generated fractional ideals of $R$
contained in $I$; and $I_w:=\bigcup(I\colon J)$ where the union is
taken over all finitely generated ideals $J$ of $R$ with
$J^{-1}=R$. We say that $I$ is {\em divisorial} if $I_v=I$; a {\em
$t-$ideal} if $I_t=I$; and a {\em $w-$ideal} if $I_w=I$. Any
divisorial ideal is a $w-$ideal. Now, $R$ is said to be a {\em
Mori domain} if it satisfies the ascending chain condition on
divisorial ideals \cite{B2,B3,BH,GH} and a {\em strong Mori
domain} if it satisfies the ascending chain condition on
$w-$ideals \cite{FM,Mi}. Trivially, a Noetherian domain is strong
Mori and a strong Mori domain is Mori.

Finally, we say that $R$ is  {\em $t-$locally Noetherian} if
$R_{M}$ is Noetherian for each maximal $t$-ideal $M$ of $R$
\cite{Kg}. Recall that strong Mori domains are $t-$locally
Noetherian \cite[Theorem 1.9]{FM}.

Throughout, we shall use $\Spec(R)$, $\Max(R)$, and $t-\Max(R)$ to denote the sets of prime
ideals, maximal ideals, and maximal $t-$ideals, respectively, of $R$.
\bigskip

We begin by providing necessary $t-$ideal-theoretic conditions for Clifford regularity.

%%%4.1%%%Old1.10%%%%%%%%%%%%%%%%%
\begin{lemma} \label{4.1}
 Let $R$ be a Clifford regular domain.Then $I_{t}\subsetneqq R$ for each nonzero proper ideal $I$ of $R$. In particular, $\Max(R)$ = $t-\Max(R)$.
\end{lemma}
%%%%%%%%%%%%%%%%%%%%%
 \begin{proof}[\bf Proof.] Deny. Then there exists a nonzero proper finitely generated ideal $I$ of $R$ such that $I_{v} = R$. So $(I:I) = I^{-1} =R$.
Hence $(I:I^{2}) = ((I:I):I)= (R:I) =I^{-1} =R$. Since ${\overline I}$ is regular in ${\mathcal S}(R)$, then $ I = I^{2}(I:I^{2}) = I^{2}$, a contradiction by \cite[Theorem 76]{K}.
\end{proof}

Next, we state our first  theorem of this section.

%%%4.2%%%Old3.10%%%%%%%%%%%%%%%%%
\begin{theorem} \label{4.2}
 Let $R$ be a $t$-locally Noetherian domain. Then $R$ is Clifford (resp., Boole) regular if and only if $R$ is stable (resp., strongly stable). Moreover, when any one condition holds, $R$ is either a field or a one-dimensional Noetherian domain.
\end{theorem}
%%%%%%%%%%%%%%%%%%%%%
 \begin{proof}[\bf Proof.] Assume that $R$ is Clifford regular. By Lemma~\ref{4.1},
 we have $\Max(R)=t-\Max(R)$. Hence $R$ is locally Noetherian. Now,
suppose that $R$ is not stable. Then there is a nonzero ideal $I$ of
$R$ such that $I(T:I)\subsetneqq T$, where $T:=(I:I)$. So there is a maximal
ideal $M$ of $R$ containing $I$ such that $(I(T:I))_{M}\subsetneqq
T_{M}\subseteq (I_{M}:I_{M})$. Set $J:= I_{M}(I_{M}:I_{M}^{2})$. By \cite[Proposition 2.9]{Ba3},
$J= (I(T:I))_{M}$. So $J$ is a nonzero proper ideal of
$(I_{M}:I_{M})$. Since ${\mathcal S}(R_{M})$ is Clifford, then $\overline{I_M}$ is regular
in ${\mathcal S}(R_{M})$. So $I_{M}=I_{M}^{2}(I_{M}:I_{M}^{2})=I_{M}J$. Since $R_{M}$
is Noetherian, then $I_{M}$ is a f.g. ideal of $R_{M}$ and therefore a f.g.
ideal of $(I_{M}:I_{M})$. By \cite[Theorem 76]{K}, $J=(I_{M}:I_{M})$, the desired contradiction.
The converse is handled by  Lemma~\ref{2.1}.

The Boolean statement follows from the Clifford statement and Lemma~\ref{3.2}.

Finally, one may assume that $R$ is a stable domain that is not a
field. Then $R$ has finite character and hence is locally
Noetherian by Lemma~\ref{4.1}. So $R$ is Noetherian by \cite[Lemma
37.3]{G}. Further, we have $\dim(R)=1$ by \cite[Proposition
2.4]{AHP}, completing the proof of the theorem.
\end{proof} Thus, a strong Mori domain that is Clifford regular
(equivalently, stable) is necessarily a Noetherian domain. Here,
Clifford regularity forces the $w-$operation to be trivial (see
also \cite[Proposition 1.3]{Mi}). Also noteworthy is that while
{\em a $t-$locally Noetherian stable domain is necessarily a
one-dimensional L-stable domain}, the converse does not hold in
general. For instance, consider an almost Dedekind domain which is
not Dedekind and appeal to Corollary~\ref{2.6}. However, the
equivalence holds for Noetherian domains:

%%%4.3%%Old3.2%%%%%%%%%%%%%%%%%%
\begin{corollary}\label{4.3} {\rm (\cite[Theorem 2.1]{Ba3} and \cite[Proposition 2.4]{AHP})}
Let $R$ be a Noetherian domain that is not a field. The following statements are equivalent:\\
i) $R$ is Clifford regular;\\
ii) $R$ is stable;\\
iii) $R$ is L-stable with $\dim(R)=1$.
\end{corollary}

%%%4.4%%Old3.5%%%%%%%%%%%%%%%%%%
\begin{corollary}\label{4.4}
Let $R$ be a local Noetherian domain such that the extension $R\subseteq {\overline R}$ is maximal, where ${\overline R}$ denotes the integral closure of $R$. The following statements are equivalent:\\
i) $R$ is Boole regular;\\
ii) $R$ is strongly stable;\\
iii) $R$ is stable and ${\overline R}$ is a PID.
\end{corollary}
%%%%%%%%%%%%%%%%%%%%%
\begin{proof}[\bf Proof.] In view of Theorem~\ref{4.2} and Proposition~\ref{2.4}, we need only prove the implication $(iii)\Longrightarrow (ii)$. Let $I$ be a nonzero ideal of $R$ and $T:=(I:I)$.  Since here ${\overline R}$ is identical to the complete integral closure of $R$, then $R\subseteq T\subseteq {\overline R}$, hence either $R=T$ or $T={\overline R}$. If $R=T$, then $I$ is invertible and hence principal in $R$ (since $R$ is local). If $T={\overline R}$, the conclusion is trivial. \end{proof}

Corollary~\ref{4.4} generates new families of Boole regular domains (i.e., with regard to those integrally closed provided by Lemma~\ref{3.5} and Corollary~\ref{3.9}).

%%%4.5%%%%%%%%%%%%%%%%%%%%
\begin{example}\label{4.5}
Let $k$ be a field and $X$ an indeterminate over $k$. \\
Let $R:=k[X^{2}, X^{3}]_{R\setminus(X^{2}, X^{3})}$. Clearly,
${\overline R}=k[X]_{R\setminus(X^{2}, X^{3})}$ is a PID and the
extension $R\subseteq {\overline R}$ is maximal. Further, $R$ is a
Noetherian Warfield domain, hence stable (cf. \cite{BS2}).
Consequently, $R$ is a one-dimensional non-integrally closed local
Noetherian domain that is Boole regular by Corollary~\ref{4.4}.
\end{example}

At this point, note that a Noetherian domain that is Clifford regular (equivalently, stable)
 need not be Boole regular (equivalently, strongly stable). For instance, consider Dedekind
  domains that are not PIDs (cf. Remark~\ref{3.7}). The following is an example of a
  non-integrally closed Noetherian Clifford regular domain that is not Boole regular.
  It also shows that Corollary~\ref{4.4} fails, in general, when $R$ is no longer local.

%%%4.6%%Old3.7%%%%%%%%%%%%%%%%%%
\begin{example}\label{4.6}
Under the same notation of the above example, let $R:=k[X^{2}, X^{3}]$. Clearly, ${\overline R}=k[X]$ and the extension $R\subseteq {\overline R}$ is maximal. Similarly,  $R$ is stable (and hence Clifford regular). However, $R$ is
not Boole regular since the ideal $I:=(X^{2}-1,X^{3}-1)$ is not principal in $(I:I)=R$.
\end{example}

We now aim toward a possible characterization of Mori domains with Clifford or Boolean class semigroup that links them to stability. In what follows, we shall use ${\overline R}$ and $R^{*}$ to denote the integral closure and complete integral closure, respectively, of an integral domain $R$. Suitable background on Mori domains is \cite{B3}. \bigskip

Next, we announce our second theorem of this section.

%%4.7%%%%%%%%%%%%%%%%%%%%%
\begin{theorem} \label{4.7}
Let $R$ be a Mori domain. The following statements are equivalent:\\
i) $R$ is a one-dimensional Clifford (resp., Boole) regular domain and $R^*$ is Mori;\\
ii) $R$ is stable (resp., strongly stable).
\end{theorem}

The proof requires the following result which provides a classification for Mori stable domains.

%%4.8%%%%%%%%%%%%%%%%%%%%%
\begin{lemma} \label{4.8}
Let $R$ be an integral domain. The following statements are equivalent:\\
i) $R$ is a Mori stable domain;\\
ii) $R$ has finite character and $R_M$ is a DVR or a one-dimensional Mori stable domain for each $M\in \Max(R)$.
\end{lemma}
%%%%%%%%%%%%%%%%%%%%
\begin{proof}[\bf Proof.] Combine \cite[Corollary 2.7]{O2} and \cite[Theorem 4.18]{GH}. \end{proof}

%%%%%%%%%%%%%%%%%%%%
\begin{proof}[\bf Proof of Theorem~\ref{4.7}.] $i)\Longrightarrow ii)$ By Proposition~\ref{2.4}, ${\overline R}$ is a Pr\"ufer domain. It follows that $R^*$ is a Dedekind domain. Further, $\dim(R)=1$ implies that $\dim_{v}(R) = 1$ by \cite[Theorem 1.10]{ABDFK}, where $\dim_{v}(R)$ denotes the valuative dimension of $R$ . Now, let
$I$ be a nonzero proper ideal of $R$. Set  $B := (I:I)$ and $J:= I(B:I)$. Suppose that $J$ is a {\em proper} ideal of $B$. Since $R\subseteq
B\subseteq R^{*}$, then $1 = \dim_{v}(R) \geq \dim_{v}(B) \geq \dim(B) \geq
1$, whence $\dim(B) = 1$. Let $P$ be a prime  ideal
of $B$ such that $J\subseteq P$. So $\htt P =1$. By \cite[Proposition 1.1]{BH}, there exists a prime ideal $Q$ of $B^{*}=R^{*}$ such that
$Q\cap B = P$. Since ${\overline I}$ is regular in ${\mathcal S}(R)$, then $I =
I^{2}(I:I^{2}) = I^{2}((I:I):I) = I^{2}(B:I) = IJ$. Hence $B = (I:I) =
(I:IJ) = ((I:I):J)= (B:J) = (J:J)$ (since $J$ is a trace ideal of $B$). So
$(J:J^{2}) = ((J:J):J) = (B:J) =B$. Hence $J^{2}(J:J^{2}) = J^{2}$. Since
$(R:B)\not =(0)$,  $B$ is Clifford regular by Lemma~\ref{2.3}. So that $J= J^{2}(J:J^{2}) = J^{2}$, hence
$JR^{*} = J^{2}R^{*}=  (JR^{*})^{2}$. Since
$R^{*}$ is a Noetherian domain,  $JR^{*} = R^{*}$ by \cite[Theorem 76]{K}, whence $R^{*} = JR^{*}\subseteq PR^{*}\subseteq Q$, absurd. Therefore $J = B$ and hence $I$ is
stable.

The Boolean statement follows from the Clifford statement and Lemma~\ref{3.2} to complete the proof of the forward direction.\\
$ii)\Longrightarrow i)$ Lemma~\ref{4.8} yields $\dim(R) =1$. It remains to show that $R^*$ is Mori, equivalently, Dedekind. Recall first that every overring of a stable domain is stable \cite[Theorem 5.1]{O3}. Thus, ${\overline R}$ is now a one-dimensional integrally closed stable domain. Hence ${\overline R}$ is Dedekind and so is $R^*$, completing the proof of the theorem.
\end{proof}

It is worth recalling that for a Noetherian domain $R$ we have: ``$\dim(R)=1$ if and only if $\dim(R^{*})=1$ if and only if $R^*$ is Dedekind'' (since here $R^{*}={\overline R}$). The same result holds if $R$ is a Mori domain such that $(R:R^{*})\not= 0$ \cite[Corollary 3.4(1) and Corollary 3.5(1)]{BH}. Also it was stated that the ``only if'' assertion holds for seminormal Mori domains \cite[Corollary 3.4(2)]{BH}. However, beyond these contexts, the problem remains elusively open. This explains the cohabitation of ``$\dim(R)=1$'' and ``$R^*$ is Mori'' hypotheses in Theorem~\ref{4.7}. In this vein, we set the following open question: ``{\tt Let $R$ be a local Mori Clifford regular domain. Is $\dim(R)=1$ if and only if  $R^*$ is Dedekind?}''\bigskip

Next, we announce our third theorem of this section. It partly draws on Theorem~\ref{4.7} and treats two well-studied large classes of Mori domains \cite{B3}. Recall that a domain $R$ is seminormal if $x\in R$ whenever $x\in K$ and $x^{2}, x^{3}\in R$ (equivalently, $x^{n}\in R$ for all $n\gg 0$).

%%4.9%%%Old3.9&3.11&3.10%%%%%%%%%%%%%%%%%%
\begin{theorem} \label{4.9}
Let $R$ be a Mori domain. Assume that either (a), (b), or (c) holds:\\
{\em a)} The conductor $(R:R^{*})\not= 0$ ; {\em b)} $R$ is seminormal ; {\em c)} The extension $R\subseteq R^{*}$ has at most one proper intermediate ring.\\
 Then $R$ is a Clifford (resp., Boole) regular domain if and only if $R$ is a stable (resp., strongly stable) domain.
\end{theorem}

The proof of (c) requires the following technical lemma.

%%4.10%%%%%%%%%%%%%%%%%%%%%
\begin{lemma} \label{4.10}
Let $R$ be a Clifford regular domain and let $I$ be a nonzero ideal of $R$. If  $(I:I)$ is a Mori domain, then $I$ is a stable ideal of $R$.
\end{lemma}
\begin{proof}[\bf Proof.]  Assume $T:=(I:I)$ is a Mori domain. By Lemma~\ref{2.3}, $T$ is Clifford regular. Suppose that $I$ is not stable. Then $J:=I(T:I)$ is a proper trace ideal of $T$. Since $\bar{I}$ is regular in ${\mathcal S}(R)$, then $I=I^{2}(I:I^{2})=IJ$. So $T=(I:I)=(I:IJ)=((I:I):J)=(T:J)=(J:J)$. Hence $J_{v}=T$. Since $T$ is Mori, then $J_{t}=J_{v}=T$ (the $v-$ and $t-$operations being with respect to $T$). Lemma~\ref{4.1} leads to the desired contradiction. \end{proof}

%%%%%%%%%%%%%%%%%%%%%%%%%%%%%%%%%%%%%%%%%%
\begin{proof}[\bf Proof of Theorem~\ref{4.9}.] We need only prove the ``only if'' assertion for Clifford regularity. Let $R$ be a Mori Clifford regular domain that is not a field. By Proposition~\ref{2.4}, $R^{*}$ is a Pr\"ufer domain.\\
a) Assume $(R:R^{*})\not= 0$. By \cite[Corollary 18]{B1}, $R^{*}$ is a Krull domain and thus Dedekind, so that $\dim(R^{*}) = 1$. It follows that $\dim(R) = 1$ by \cite[Corollary 3.4]{BH}. Theorem~\ref{4.7} leads to the conclusion.\\
b) Assume that $R$ is seminormal. According to \cite[Theorem 2.9]{B2}, $R^{*}$ is a Krull domain and hence Dedekind. In view of Theorem~\ref{4.7}, we need only show that $\dim(R)=1$. Let $M$ be any  maximal ideal of $R$. Clearly, $R_M$ is a seminormal local Mori Clifford regular domain. Therefore, we may assume that $R$ is local with maximal ideal $M$. Suppose that $\htt M=\dim(R)\geq
2$. By Lemma~\ref{4.1}, $M$ is a $t-$ideal of $R$. Since
$R$ is Mori, then $M_{v}=M_{t}=M$. Hence $R\subsetneqq M^{-1}$. By \cite[Proposition 1]{Q}, $M$ is strongly divisorial, so that $T:= (M:M)=M^{-1}$ is a Mori
domain. Since $R$ is seminormal, by \cite[Lemma 2.5]{BH} there is a non-divisorial prime $Q$ of $T$ contracting on $M$ such that $\htt _{T}Q\geq 2$. Since $Q$ is not divisorial in $T$,
$(T:Q)=T$  by \cite[Proposition 1]{Q}, whence $Q_{t}=Q_{v}=T$ (the $t-$ and $v-$operation being with respect to $T$). Further, $T$ is Clifford regular by Lemma~\ref{2.3}. Therefore, Lemma~\ref{4.1} yields the desired contradiction. Hence $\dim(R)=1$, as desired.\\
c) Assume that $R\subseteq R^{*}$ has at most one proper intermediate ring. Let $I$ be a nonzero ideal of $R$ and let $J:=II^{-1}$. Since $R\subseteq
(I:I)\subseteq R^{*}$, then either $(I:I)=R^{*}$, $R=(I:I)$, or $R\subsetneqq (I:I)\subsetneqq R^{*}$. In view of (a) and Lemma~\ref{4.10}, we need only handle the late case. Since now
$R\subsetneqq (I:I)\subseteq (J:J)=J^{-1}\subseteq R^{*}$, then either $(J:J)=J^{-1}=R^{*}$ or $(I:I)=(J:J)=J^{-1}$. The former case follows from (a). The latter case follows from Lemma~\ref{4.10}, since $J^{-1}$ is a Mori domain by \cite[Theorem 11]{Lu}. Consequently, in all cases $I$ is stable and so is $R$.  \end{proof}

One may wonder about the existence of (one-dimensional) Mori stable domains that are not Noetherian. Indeed, the pullback construction -a main source for non-Noetherian non-Krull Mori domains- can be of no help in this regard. More precisely, let $T$ be a domain, $M$ a maximal ideal of $T$, $K$ its residue field, $\phi:T\longrightarrow K$ the canonical surjection, and $D$ a proper subring of $K$ with quotient field $\qf(D)=k$. Let $R:=\phi^{-1}(D)$. Then $R$ is a Mori stable domain only if $R=T$.
This follows easily from a combination of \cite[Theorem 4.18]{GH} and \cite[Theorem 2.6]{O2}.
(i.e., while the former result yields $D=k$, the latter, applied to ($T_M$, $R_M$, $MR_M$), yields $k=K$.)

Also, it turns out that non-Noetherian Mori Clifford regular domains can't stem
from our results on pullbacks (Section 5). Indeed, under the hypotheses of
Theorem~\ref{5.1}(2) below, Noetherianity and the Mori property coincide for the pullback $R$.

%%%%%%%%%%%%%%%%%%%%%%%%%%%%%%%%%%%%%%%%%%%%%%%%%%%%%%%%
%%%SECTION5%%%%%%%%%%%%%%%%%%%%%%%%%%%%%%%%%%%%%%%%%%%%%
%%%%%%%%%%%%%%%%%%%%%%%%%%%%%%%%%%%%%%%%%%%%%%%%%%%%%%%%
\section{Pullbacks}

The purpose of this section is to characterize Clifford regularity and Boole
regularity in pullback constructions. Our work is motivated by an attempt to
generating  new families of integral domains with Clifford or Boolean class
semigroup, moving therefore beyond the classical contexts of integrally
closed or Noetherian domains.\bigskip

Let us fix the notation for the rest of this section. Let $T$ be an integral
domain, $M$ a maximal ideal of $T$, $K$ its residue field,
$\phi:T\longrightarrow K$ the canonical surjection, $D$ a proper subring of
$K$, and $k:=\qf(D)$. Let $R:=\phi^{-1}(D)$ be the pullback issued from the
following diagram of canonical homomorphisms:
\[\begin{array}{ccl}
R            & \longrightarrow                 & D\\
\downarrow   &                                 & \downarrow\\
T            & \stackrel{\phi}\longrightarrow  & K=T/M
\end{array}\]\medskip

First, we wish to shed light on some features imposed by a possible passage
of Clifford regularity to pullbacks. As a matter of fact, $R$ need not be
Clifford regular even when $D$ is a PID with $k=K$ and $T$ is a DVR (e.g.,
$R:=\Z+X\Q[[X]]$) or  when $D=k$ and $T$ is local (see Example~\ref{5.2}).
In the well-studied case where $T$ is integrally closed (e.g., a valuation
domain or a polynomial ring over a field), Clifford regularity of $R$
transfers to ${\overline R}$, since  here ${\overline R}=\phi^{-1}(D')$,
where $D'$ is the integral closure of $D$ in $K$. Further, ${\overline R}$
and (hence) $R$ have finite character, which forces $D$ to be semilocal.
This follows easily from a combination of Lemma~\ref{2.3},
Proposition~\ref{2.4}, the finite character requirement \cite[Theorem
4.5]{Ba3}, and the well-known fact that $\Spec(R)$ is an amalgamated sum of
$\Spec(D)$ and $\Spec(T)$ over the conductor $M$ \cite{F}. \bigskip

Next, we announce our first theorem of this section. It particularly
provides a necessary and sufficient condition for a pseudo-valuation domain
(i.e., PVD) to inherit Clifford or Boole regularity.

%%5.1%%old4.1%%%%%%%%%%
\begin{theorem} \label{5.1} Under the above notation, the following hold:\\
1) If $R$ is Clifford (resp., Boole) regular, then so are $T$ and $D$, and
$[K:k]\leq 2$. \\
2) Assume $D=k$ and $T$ is a valuation (resp., strongly discrete valuation)
domain. Then $R$ is Clifford (resp., Boole) regular if and only if
$[K:k]=2$.
\end{theorem}
%%%%%%%%%%%%%%%%%%%%%%%%

We need the following technical lemma.

%%5.2%%%%%%%%%%%%%%%%%%%%%
\begin{lemma} \label{5.1(2)} Under the above notation, let $W$ be a $D-$submodule of $K$ containing $D$. Then $\phi^{-1}(W:W)=(\phi^{-1}(W):\phi^{-1}(W))$.
\end{lemma}
%%%%%%%%%%%%%%%%%%%%%%%%
\begin{proof}[\bf Proof.] Let $W$ be a $D$-module such that $D\subseteq W\subsetneqq
K$. Since $1\in W$, then $(W:W)\subseteq W$. So $\phi^{-1}(W:W)\subseteq
\phi^{-1}(W)\subseteq T$. Now, let $x\in \phi^{-1}(W:W)$. So, for each $z\in \phi^{-1}(W)$, $\phi(xz)=\phi(x)\phi(z)\in W$. Then $xz\in \phi^{-1}(W)$ and therefore
$x\in (\phi^{-1}(W):\phi^{-1}(W))$. Conversely, let $x\in
(\phi^{-1}(W):\phi^{-1}(W))$. Since $1\in \phi^{-1}(W)$, then $x\in
\phi^{-1}(W)\subseteq T$ and $x\phi^{-1}(W)\subseteq \phi^{-1}(W)$ implies
that $\phi(x)W=\phi(x\phi^{-1}(W))\subseteq \phi(\phi^{-1}(W))=W$. Hence
$\phi(x)\in (W:W)$, as desired. \end{proof}

%%%%%%%%%%%%%%%%%%%%%%%%%%%%%%%%%%%%
\begin{proof}[\bf Proof of Theorem~\ref{5.1}.] 1) Assume that $R$ is Clifford (resp., Boole) regular. Then so is $T$ by Lemma~\ref{2.3}. Let $J$ be a nonzero
(integral) ideal of $D$ and let $I:=\phi^{-1}(J)$. By \cite[Proposition 6]{HKLM}, $(I:I^{2})=\phi^{-1}(J:J^{2})$. So
$J=\phi(I)=\phi(I^{2}(I:I^{2}))=J^{2}(J:J^{2})$ and therefore $D$ is
Clifford regular. Now, assume that $R$ is Boole regular. Then there exists
$0\not =c\in \qf(R)$ such that $I^{2}=cI$. Since $J$ is nonzero, then
$M\subsetneqq I$. Let $R_{0}=\phi^{-1}(k)$  be the pullback issued from the
following diagram of canonical homomorphisms:
\[\begin{array}{ccl}
R            & \longrightarrow                 & D\\
\downarrow   &                                 & \downarrow\\
R_{0}& \longrightarrow                 & k\\
\downarrow   &                                 & \downarrow\\
T            & \stackrel{\phi}\longrightarrow  & K=T/M
\end{array}\]
Since $M\subsetneqq I\subseteq IR_{0}$ and $M$ is a maximal ideal of
$R_{0}$, then $IR_{0}=R_{0}$. So
$1=\displaystyle\sum_{i=1}^{i=n}a_{i}x_{i}$, where $a_{i}\in I$ and
$x_{i}\in R_{0}$ for each $i$. Hence
$c=\displaystyle\sum_{i=1}^{i=n}ca_{i}x_{i}$. Since $ca_{i}\in
cI=I^{2}\subseteq R\subseteq R_{0}$, then $ca_{i}x_{i}\in R_{0}$ for each
$i$, hence $c\in R_{0}$. So $\phi(c)\in k=\qf(D)$ and
$J^{2}=\phi(c)J$. It follows that $D$ is Boole regular. It remains to prove that $[K:k]\leq 2$. Notice first that $R_{0}$ is Clifford by Lemma~\ref{2.3}. \\
{\bf Step 1}. We claim that, for each $x\in K$, $x^{2}\in k+xk$. By
a contrast way, suppose there exists $x\in K$ such that $x^{2}\not\in k+xk$. let
$W$ be the $k$-vector space defined by $W:=k+xk$ and let $I$ be the ideal of $R_{0}$ given by $I :=m\phi^{-1}(W)$ for some nonzero $m\in M$. We first show that $(W:W) =k$.
It is clear that $k\subseteq (W:W)$. Since $1\in W$, then $(W:W)\subseteq
W$. Let $z\in (W:W)$. Write $z= a+bx$, where $a, b\in k$. Since
$x\in W$, then $zx\in W$. So $bx^{2}+ax=zx=c+dx$ for some $c, d\in k$. If
$b\not =0$, then $x^{2}=b^{-1}(d-a)x + b^{-1}c\in k+xk$, which is absurd. So
$b=0$ and therefore $z=a\in k$. Hence $(W:W)=k$. Now, by Lemma~\ref{5.1(2)},
$(I:I)=(m\phi^{-1}(W):m\phi^{-1}(W))=
(\phi^{-1}(W):\phi^{-1}(W))=\phi^{-1}((W:W))=\phi^{-1}(k)=R_{0}$. So
$(I:I^{2})=((I:I):I)=(R_{0}:I)=m^{-1}\phi^{-1}((k:W))=m^{-1}\phi^{-1}(0)=m^{-1}M$.
Hence $I^{2}(I:I^{2})\subseteq mM\subsetneqq I$, which is a contradiction since
${\overline I}$ is regular in ${\mathcal S}(R_{0})$. It follows that for each $x\in
K\setminus k$, $[k(x):k]=2$.\\
{\bf Step 2}. Suppose that $[K:k]\geq 3$. Consider a free system
$\{1, x, z \}$ of $K$ as a $k$-vector space. Let $W:=k+xk+zk$ and
$I:=m\phi^{-1}(W)$ for some nonzero $m\in M$. We whish to show that $(W:W)=k$. Let $y\in (W:W)\subseteq
W$. Write $y= a+bx+cz$. Since $x\in W$, then $xy= ax +bx^{2}+czx\in W$. By
the first step, $x^{2}= dx+e$ for some $d, e\in k$. Hence
$ax+bdx+be +cxz=xy\in W$. So $cxz=xy-(a+bd)x-be\in W$. If $c\not =0$, then
$xz\in W$, whence $W$ is a ring. So $W=k[x,z]=k(x,z)$ (since, by the first
step, $x$ and $z$ are algebraic over $k$). Hence $[W:k]=[k(x, z):k] =4$
which is absurd. It follows that $c=0$. Similarly, using the fact that
$z\in W$, we obtain that $b=0$. Hence $y=a\in k$ and therefore $(W:W) =k$.
Now, as in the first step, we obtain that $I^{2}(I:I^{2})\subseteq mM\subsetneqq
I$, which is a contradiction. It follows that $[K:k]=2$.\\
2) Assume that $D=k$ and $[K:k]=2$. Let $I$ be a nonzero (integral) ideal of
$R$. If $I$ is an ideal of $T$, since $T$ is Clifford (resp., Boole)
regular, then $I^{2}(I:I^{2})=I$ (resp., $I^{2}=cI$). If $I$ is not an ideal
of $T$, then as in \cite[Theorem 1]{BG}, it is easy to see that $I=c\phi^{-1}(W)$,
where $k\subseteq W\subsetneqq K$ is a $k$-vector space. Since $[K:k]=2$, then
$W=k$ and therefore $I=cR$, as desired. \end{proof}

The following example shows that  Theorem~\ref{5.1}(2) does not hold in
general, and hence nor does the converse of (1).

%%5.3%%old4.2%%%%%%%%%%%%%%%%%%%
\begin{example} \label{5.2}
Let $\Z$ and $\Q$ denote the ring of integers and field of rational numbers,
respectively, and let $X$ and $Y$ be indeterminates over $\Q$.
Set  $V:=\Q(\sqrt{2},\sqrt{3})[[X]]$, $M:=X\Q(\sqrt{2},\sqrt{3})[[X]]$,
$T:=\Q(\sqrt{2})+M$, and $R:=\Q+M$. Both $T$ and $R$ are one-dimensional
local Noetherian domains arising from the DVR $V$, with ${\overline T}=V$
and ${\overline R}=T$. By Theorem~\ref{5.1}(2), $T$ is Clifford (actually,
Boole) regular, whereas $R$ is not. More specifically, the isomorphy class
of the ideal $I:=X(\Q+\sqrt{2}\Q+\sqrt{3}\Q+M)$ is not regular in ${\mathcal
S}(R)$.
\end{example}

The following is an immediate consequence of Theorem~\ref{5.1}(2).

%%5.4%%%%%%%%%%%%%%%%%%%%%
\begin{example} \label{5.3} Let $n$ be an integer $\geq 1$. Let $R$ be a PVD
associated to a non-Noetherian valuation (resp., strongly discrete
valuation) domain $(V,M)$ with $\dim(V)=n$ and $[V/M:R/M]=2$. Then $R$ is an
$n-$dimensional  local Clifford (resp., Boole) regular domain that is
neither integrally closed nor Noetherian.
\end{example}

Next, we provide new examples of Noetherian Boole (hence Clifford) regular
domains (with regard to Example~\ref{4.5}).

%%5.5%%%%%%%%%%%%%%%%%%%%%
\begin{example} \label{5.4}  Let $R$ be a PVD associated to a DVR $(V,M)$
with $[V/M:R/M]=2$. Then $R$ is a one-dimensional local Noetherian Clifford
Boole regular domain that is not integrally closed.
\end{example}

Now, we introduce a useful class of domains that may help constructing more
original examples for Clifford or Boole regularity. An integral domain $A$
is said to be {\em conducive} if the conductor $(A:B)$ is nonzero for each
overring $B$ of $A$ other than its quotient field. Examples of conducive
domains include pseudo-valuation domains and, in general, arbitrary
pullbacks of the form $R:=D+M$ arising from a valuation domain $V:=K+M$
\cite[Propositions 2.1 \& 2.2]{DF}. Suitable background on conducive domains
is \cite{BDF,DF}. \bigskip

We are now able to announce our second theorem of this section. It
treats Clifford regularity, for the remaining case ``$k=K$'', for
pullbacks $R:=\phi^{-1}(D)$ where $D$ is a conducive domain.

%%5.6%%old4.5and 4.8%%%%%%%%%%%%%%%%%%%
\begin{theorem} \label{5.5} Under the above notation, suppose that $D$ is a
semilocal conducive domain with quotient field $k=K$ and  either (a) or (b)
holds:\\
a) $T$ is a valuation domain, \\
b) $T:=K[X]$ and $R:=D+XK[X]$, where $X$ is an indeterminate over $K$.\\
Then $R$ is Clifford regular if and only if so is $D$.
\end{theorem}

The proof of (a) is actually handled by the following technical lemma.

%%5.7%%old4.4%%%%%%%%%%%%%%%%%%%
\begin{lemma} \label{5.6} Under the above notation, suppose that $T$ is a
valuation domain and for each $D-$submodule $W$ of $K$ containing $D$,
either $W$ is a ring or $(D:W)\not=0$. \\
Then $R$ is Clifford regular if and only if so is $D$.
\end{lemma}
%%%%%%%%%%%%%%%%%%%%%%%%
\begin{proof}[\bf Proof.] We need only prove the ``if'' assertions. Assume that $D$ is Clifford regular. Let $I$ be a nonzero (integral) ideal of $R$. If $M\subsetneqq I$, then
$I=\phi^{-1}(J)$ for some nonzero ideal $J$ of $D$. Since $D$ is Clifford
regular, then $J^{2}(J:J^{2})=J$. By \cite[Proposition 6]{HKLM}, it is easy to see
that $I^{2}(I:I^{2})=I$. Assume that $I\subseteq M$. If $I$ is an ideal of
$T$, we are done (since $T$ is Clifford regular). If $I$ is not an
ideal of $T$, then as in \cite[Theorem 1]{BG}, it is easy to see that
$I=c\phi^{-1}(W)$, where $W$ is a $D$-module with $D\subseteq W\subsetneqq K$.
If $W$ is a ring, then clearly $W^{2}(W:W^{2})=W$ and therefore
$I^{2}(I:I^{2})=I$ by Lemma~\ref{5.1(2)}. If $(D:W)\not =(0)$, then $dW$ is an
(integral) ideal of $D$ for some nonzero element $d$ of $D$. Since $D$ is
Clifford regular, then $(dW)^{2}(dW:(dW)^{2})=dW$ so that $W^{2}(W:W^{2})=W$. Therefore $I^{2}(I:I^{2})=I$ by Lemma~\ref{5.1(2)}.
\end{proof}

%%%%%%%%%%%%%%%%%%%%%%%%
\begin{proof}[\bf Proof of Theorem~\ref{5.5}.] a) Follows easily from Lemma~\ref{5.6}.\\
b) Assume that $D$ is Clifford regular. Let $I$ be a nonzero ideal of $R$.
Then $I=f(X)(F+XK[X])$, where $F$ is a nonzero $D$-submodule of $K$ such
that $f(0)F\subseteq D$ \cite[Proposition 4.12]{CMZ}. Since $D$ is conducive,
then $F$ is a fractional ideal of $R$. Hence $F^{2}(F:F^{2})=F$ and
therefore $I^{2}(I:I^{2})=I$, as desired.
\end{proof}

Clearly, Theorems \ref{5.1} and \ref{5.5} generate new families of examples
of Clifford regular domains, as shown by the following construction.

%%5.8%%%%%%%%%%%%%%%%%%%%%
\begin{example} \label{5.7} For every positive integer $n\geq 2$, there
exists an example of an integral domain $R$ satisfying the following
conditions:\\
1) $\dim(R)=n$,\\
2) $R$ is neither integrally closed nor Noetherian,\\
3) $R$ is Clifford regular,\\
4) Each overring of $R$ is Clifford regular,\\
5) $R$ has infinitely many maximal ideals.
\end{example}
%%%%%%%%%%%%%%%%%%%%%%%%
\begin{proof}[\bf Proof.] Here is an explicit example. Let $n\geq 1$ and let $X,X_{1}, ...,
X_{n-1}$ be indeterminates over $\Q$. Set
$V_{1}:=\Q(\sqrt{2})+M_{1}$, where
$M_{1}:=X_{1}\Q(\sqrt{2})[X_{1}]_{(X_{1})}$;
$V_{i}:=V_{i-1}+M_{i}$, where $M_{i}:=X_{i}\Q(\sqrt{2},X_{1}, ...,
X_{i-1})[X_{i}]_{(X_{i})}$ for each $2\leq i\leq n-1$; $M:=M_{1} +
... + M_{n-1}$; $D:=\Q+M$; $K=\Q(\sqrt{2},X_{1}, ..., X_{n-1})$;
and $R:=D+XK[X]$. Clearly, $V:=V_{n-1}=\Q(\sqrt{2})+M$ is an
$(n-1)-$dimensional valuation domain with maximal ideal $M$
\cite[Theorem 2.1]{BG}, ${\overline R}:=V+XK[X]$, and hence $R$ is
an $n-$dimensional non-integrally closed non-Noetherian domain
\cite{AD,BG,BR,CMZ,G}. Further, $R$ is Clifford regular by
Theorem~\ref{5.1} and Theorem~\ref{5.5}. Now let $S$ be an
overring of $R$. Since $V\subseteq{\overline S}$ and
$\qf(D)=\qf(V)=K$, it easily can be seen that $V\subseteq S$,
hence ${\overline R}\subseteq S$. Consequently, $S$ is Clifford
regular since ${\overline R}$ is. Finally, $\Spec(R)$ has the
following shape \cite{AD,BG,CMZ}:

%\newpage
%%DIAGRAM%%%%%%%%%%%%%%%%%%%%%%%%%%
%%%%%%%%%%%%%%%%%%%%%%%%%%%%%%%%%%%
\[\setlength{\unitlength}{1mm}
\begin{picture}(80,60)(0,10)
%%%%%%%%%%%%
\put(40,70){\line(0,-1){10}}
\put(40,60){\line(0,-1){10}}
\put(40,40){\line(0,-1){10}}
\put(40,30){\line(0,-1){10}}
\put(40,20){\line(0,-1){10}}
\put(40,10){\line(1,2){5}}
\put(40,10){\line(1,1){10}}
\put(40,10){\line(3,2){15}}
%%%%%%%%%%%%%%%%%
\put(40,70){\circle*{.7}}
\put(38,70){\makebox(0,0)[r]{$XK[X]+M$}}
\put(40,60){\circle*{.7}}
\put(38,60){\makebox(0,0)[r]{$XK[X]+M_{1} + ... + M_{n-2}$}}
\put(40,50){\circle*{.7}}
\put(40,45){\circle*{.7}}
\put(40,40){\circle*{.7}}
\put(40,30){\circle*{.7}}
\put(38,30){\makebox(0,0)[r]{$XK[X]+M_{1}$}}
\put(40,20){\circle*{.7}}
\put(38,20){\makebox(0,0)[r]{$XK[X]$}}
\put(45,20){\circle*{.7}}
\put(50,20){\circle*{.7}}
\put(55,20){\circle*{.7}}
\put(60,20){\circle*{.7}}
\put(65,20){\circle*{.7}}
\put(70,20){\circle*{.7}}
\put(40,10){\circle*{.7}}
\put(38,10){\makebox(0,0)[r]{(0)}}
\end{picture}\]\end{proof}
%%%%%%%%%%%%%%%%%%%%%%%%%%%%%%%%%%%%%%%%%%%%%%%%%%%%%%%%%

%%%REFERENCES%%%%%%%%%%%%%%%%%%%%%%%%%%%%%%%%%%%%%%%%%%%%
%%%%%%%%%%%%%%%%%%%%%%%%%%%%%%%%%%%%%%%%%%%%%%%%%%%%%%%%%

%%%%%%%%%%%%%%%%%%%%%%%%%%%%%%%%%%%%%%%%%%%%%%%%%%%%%%%%%
\end{document}